\documentclass{amsart}
\usepackage{amsmath}
\usepackage{array}
\usepackage{amsfonts}
\usepackage{latexsym}
\usepackage{epsfig}
\newtheorem{theorem}{Theorem}
\newtheorem{prop}[theorem]{Proposition}

\newtheorem{lemma}[theorem]{Lemma}
\newtheorem{corollary}[theorem]{Corollary}

\newcommand\beq{\begin{equation}}
\newcommand\eeq{\end{equation}}
\newcommand\bce{\begin{center}}
\newcommand\ece{\end{center}}
\newcommand\bea{\begin{eqnarray}}
\newcommand\eea{\end{eqnarray}}
\newcommand\ben{\begin{enumerate}}
\newcommand\een{\end{enumerate}}
\newcommand\bsb{\begin{Sb}}
\newcommand\esb{\end{Sb}}

\newcommand\ms{\medskip}

\newcommand\brr{\begin{array}}
\newcommand\err{\end{array}}
\newcommand\bt{\begin{tabular}}
\newcommand\et{\end{tabular}}

\renewcommand\S{{\mathcal S}}
\newcommand\D{{\mathcal D}}

\newcommand\knu{\Psi}
\newcommand\krar{\Phi}
\newcommand\fp{f}
\newcommand\exc{e}
\newcommand\lt{\alpha}
\newcommand\rt{\beta}
\newcommand\ct{\gamma}
\newcommand\he{\nu}
\newcommand\rk{r}

\newcommand\XX{\Theta}

\newenvironment{abstrac}{%
         \small
        \begin{center}%
          {\bfseries {Abstract}\vspace{-.5em}}%
                 \end{center}%
        \quotation}

\title{Bijections for refined restricted permutations}
\author{Sergi Elizalde}
\address{Department
of Mathematics, MIT, Cambridge MA
02139.}\email{sergi@math.mit.edu}
\author{Igor Pak}
\address{
Department
of Mathematics, MIT, Cambridge MA
02139.}\email{pak@math.mit.edu}
\date{December 20, 2002}
\begin{document}
 \maketitle
 \thispagestyle{empty}

\begin{abstrac}
\medskip
We present a bijection between 321- and 132-avoiding permutations
that preserves the number of fixed points and the number of
excedances. This gives a simple combinatorial proof of recent
results of Robertson, Saracino and Zeilberger~\cite{RSZ02}, and
the first author~\cite{Eli02}. We also show that our bijection
preserves additional statistics, which extends the previous
results.

\end{abstrac}

\section{Introduction}

The subject of {\it pattern avoiding permutations}, also called
{\it restricted permutations}, has blossomed in the past decade. A
number of enumerative results have been proved, new bijections found,
and connections to other fields established.  Despite recent
progress, the so called Stanley-Wilf conjecture giving an
exponential upper bound on the number of pattern avoiding
permutations remains open, and much of the ongoing research is
related to the conjecture.

An unexpected recent result of Robertson, Saracino and
Zeilberger~\cite{RSZ02} gives a new and exciting extension to what
is now regarded as classical result that the number of 321-avoiding
permutations equals the number of 132-avoiding permutations.  They
show that one can ``refine'' this result by taking into account
the number of fixed points in a permutation.  In fact, they study
all 6 patterns in $S_3$ which produce different ``refined''
statistics, with the above mentioned result having a
highly nontrivial and technically involved proof. The story
continued in a recent paper of the first author~\cite{Eli02} where
an additional statistic, ``the number of excedances'', was added.
The proof uses some nontrivial generating function machinery and
is also quite involved.

In this paper we present a bijective proof of the ``refined''
results on 321- and 132-avoiding permutations, resolving the problem
which was left open in~\cite{RSZ02,Eli02}.  In fact, our bijection is
a composition of two (slightly modified) known bijections into
Dyck paths, and the result follows from a new analysis of these
bijections.  The Robinson-Schensted-Knuth correspondence is a part
of one of them, and the difficulty of the analysis stems from
the complexity of this
celebrated correspondence. As a new application of our bijections,
we show that the length of the longest increasing subsequence in
321-avoiding permutations corresponds to a rank in 132-avoiding
permutations, which further refines the previous results.

\smallskip

Let $n$, $m$ be two positive integers with $m\le n$, and let
$\sigma=(\sigma(1),\sigma(2),\ldots,\sigma(n))\in\S_n$ and
$\pi=(\pi(1),\pi(2),\ldots,\pi(m))\in\S_m$. We say that $\sigma$
\emph{contains} $\pi$ if there exist indices $i_1<i_2<\ldots<i_m$
such that $(\sigma(i_1),\sigma(i_2),\ldots,\sigma(i_m))$ is in the
same relative order as $(\pi(1),\pi(2),\ldots,\pi(m))$. If
$\sigma$ does not contain $\pi$, we say that $\sigma$ is
\emph{$\pi$-avoiding}. For example, if $\pi=132$, then
$\sigma=(2,4,5,3,1)$ contains~$132$, because the subsequence
$(\sigma(1),\sigma(3),\sigma(4))=(2,5,3)$ has the same relative
order as~$(1,3,2)$. However, $\sigma=(4,2,3,5,1)$ is
$132$-avoiding.

We say that $i$ is a \emph{fixed point} of a permutation $\sigma$
if $\sigma(i)=i$.  Similarly,~$i$ is an \emph{excedance} of $\sigma$
if $\sigma(i)>i$. Denote by $\fp(\sigma)$ and $\exc(\sigma)$ the
number of fixed points and the number of excedances of~$\sigma$,
respectively.

Denote by $\S_n(\pi)$ the set of $\pi$-avoiding permutations
in $\S_n$. For the case of patterns of length 3, it is
known~\cite{Knu73} that regardless of the pattern $\pi\in\S_3$,
$|\S_n(\pi)|=C_n=\frac{1}{n+1}{2n\choose n}$, the $n$-th
Catalan number. While equalities
$|\S_n(132)|=|\S_n(231)|=|\S_n(312)|=|\S_n(213)|$
and $|\S_n(321)|=|\S_n(123)|$ are straightforward, the equality
$|\S_n(321)|=|\S_n(132)|$ is more difficult to establish.
Bijective proofs of this fact are given in
\cite{Kra01,Ric88,SS85,Wes95}. However, none of these bijections
preserves either of the statistics $\fp$ or $\exc$.

\medskip

\begin{theorem}\cite{RSZ02, Eli02} \label{th:old}
\, The number of 321-avoiding permutations
$\sigma \in \S_n$ with $\fp(\sigma)=i$ and $\exc(\sigma)=j$
equals the number of 132-avoiding permutations $\sigma \in \S_n$
with $\fp(\sigma)=i$ and $\exc(\sigma)=j$, for any  $0 \le i,j\le n$.
\end{theorem}

A special case of the theorem, which ignores the number of
excedances, was given in~\cite{RSZ02}.  In full, the theorem was
shown in~\cite{Eli02}. As we mentioned above, both proofs are
non-bijective and technically involved. The main result of this
paper is a bijective proof of the following extension of
Theorem~\ref{th:old}.

\smallskip

Let $\ell(\sigma)$ be the length of the
\emph{longest increasing subsequence} of $\sigma$, i.e., the
largest $m$ for which there exist indices $i_1<i_2<\ldots<i_m$
such that $\sigma(i_1)<\sigma(i_2)<\ldots<\sigma(i_m)$. Define the
\emph{rank} of $\sigma$, denotes $\rk(\sigma)$, to be the
largest $k$ such that $\sigma(i)>k$ for all $i\le k$.

\begin{theorem}\label{th:new}
The number of 321-avoiding permutations $\sigma \in \S_n$
with $\fp(\sigma)=i$, $\exc(\sigma)=j$ and $\ell(\sigma)=k$ equals
the number of 132-avoiding permutations $\sigma \in \S_n$
with $\fp(\sigma)=i$, $\exc(\sigma)=j$ and $\rk(\sigma)=n-k$, for
any $0 \le i,j,k\le n$.
\end{theorem}

To prove this theorem, we establish a bijection~$\Theta$ between
$\S_n(321)$ and $\S_n(132)$, which respects the statistics as above.
While~$\Theta$ is not hard to define, its analysis is less
straightforward and will occupy much of the paper.

The rest of the paper is structured as follows.
In section~\ref{sec:dyck}  we define Dyck paths and several
new statistics on them. The description of the main bijection is
done in section~\ref{sec:bijection}, and is divided into two parts.
First we give a bijection from 321-avoiding permutations to
Dyck paths, and then another one from Dyck paths to 132-avoiding
permutations. In sections~\ref{sec:pf132} and \ref{sec:pf321} we
establish properties of these bijections which imply Theorem~2.
Section~\ref{sec:pflemmas} contains proofs of two technical
lemmas.  We conclude with final remarks, extensions
and open problems.

\section{Statistics on Dyck paths}\label{sec:dyck}

Recall that a \emph{Dyck path} of length $2n$ is a lattice path in
$\mathbb{Z}^2$ between $(0,0)$ and $(2n,0)$ consisting of up-steps
$(1,1)$ and down-steps $(1,-1)$ which never goes below the
$x$-axis. Sometimes it will be convenient to encode each up-step
by a letter $u$ and each down-step by $d$, obtaining an encoding
of the Dyck path as a \emph{Dyck word}. We shall denote by $\D_n$
the set of Dyck paths of length $2n$, and by
$\D=\bigcup_{n\geq0}\D_n$ the class of all Dyck paths.

For any $D\in\D$, we define a \emph{tunnel} of $D$ to be a
horizontal segment between two lattice points of $D$ that
intersects $D$ only in these two points, and stays always below
$D$. Tunnels are in obvious one-to-one correspondence with
decompositions of the Dyck word $D=AuBdC$, where $B\in\D$ (no
restrictions on $A$ and $C$). In the decomposition, the tunnel is
the segment that goes from the beginning of $u$ to the end of $d$.
If $D\in\D_n$, then $D$ has exactly $n$ tunnels, since such a
decomposition can be given for each up-step of $D$.

A tunnel of $D\in\D_n$ is called a \emph{centered tunnel} if the
$x$-coordinate of its midpoint (as a segment) is $n$, that is, the
tunnel is centered with respect to the vertical line through the
middle of $D$. In terms of the decomposition $D=AuBdC$, this is
equivalent to saying that $A$ and $C$ have the same length. Denote
by $\ct(D)$ the number of centered tunnels of $D$.

\begin{figure}[hbt]
\epsfig{file=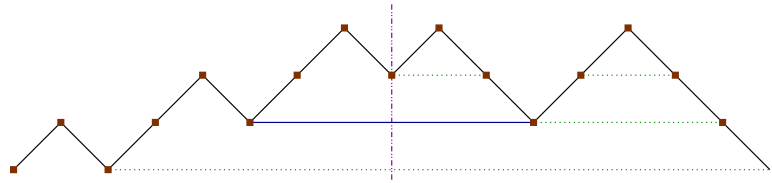,width=4in} \caption{\label{fig:ctrt}
One centered and three right tunnels.}
\end{figure}

A tunnel of $D\in\D_n$ is called a \emph{right tunnel} if the
$x$-coordinate of its midpoint is strictly greater than $n$, that
is, the midpoint of the tunnel is to the right of the vertical
line through the middle of $D$. In terms of the decomposition
$D=AuBdC$, this is equivalent to saying that the length of $A$ is
strictly bigger than the length of $C$. Denote by $\rt(D)$ the
number of right tunnels of $D$. In Figure~\ref{fig:ctrt}, there is
one centered tunnel drawn with a solid line, and four right
tunnels drawn with dotted lines. Similarly, a tunnel is called a
\emph{left tunnel} if the $x$-coordinate of its midpoint is
strictly less than $n$. Denote by $\lt(D)$ the number of left
tunnels of $D$. Clearly, $\lt(D)+\rt(D)+\ct(D)=n$ for any~$D\in\D_n$.

We will distinguish between right tunnels of $D\in\D_n$ that are
entirely contained in the half plane $x\ge n$ and those that cross
the line $x=n$. These will be called \emph{right-side
tunnels} and  \emph{right-across tunnels}, respectively.  In
Figure~\ref{fig:ctrt} there are three right-side tunnels and one
right-across tunnel. \emph{Left-side tunnels} and
\emph{left-across tunnels} are defined analogously.

Finally, for any $D\in\D_n$, define $\he(D)$ to be the height of the
middle point of $D$, that is, the $y$-coordinate of the intersection
of the path with $x=n$.

\ms

We say that $i$ is an \emph{antiexcedance} of $\sigma$ if
$\sigma(i)<i$. Sometimes it will be convenient to represent a
permutation $\sigma\in\S_n$ as an $n\times n$ array with a cross
on the squares $(i,\sigma(i))$. Note that fixed points,
excedances, and antiexcedances correspond respectively to crosses
on, strictly to the right, and strictly to the left of the main
diagonal of the array.

\section{Two bijections into Dyck paths}\label{sec:bijection}

The bijection $\Theta: \S_n(321) \longrightarrow \S_n(132)$
that we present will be the composition of two
bijections, one from $\S_n(321)$ to $\D_n$, and another one from
$\D_n$ to $\S_n(132)$.

\ms

The first bijection $\knu:\S_n(321)\longrightarrow\D_n$ is
essentially due to Knuth~\cite{Knu73}. Its definition consists of
two steps.  Given $\sigma\in\S_n(321)$, we start by applying the
Robinson-Schensted-Knuth correspondence to $\sigma$
(see e.g.~\cite{EC}).  This
correspondence gives a bijection between
the symmetric group $\S_n$ and pairs $(P,Q)$ of \emph{standard
Young tableaux} of the same shape $\lambda\vdash n$. For
$\sigma\in\S_n(321)$ the algorithm is particularly easy because in
this case the tableaux $P$ and $Q$ have at most two rows. The
\emph{insertion tableau} $P$ is obtained by reading $\sigma$ from
left to right and, at each step, inserting $\sigma(i)$ to the
partial tableau obtained so far. Assume that
$\sigma(1),\ldots,\sigma(i-1)$ have already been inserted. If
$\sigma(i)$ is larger than all the elements on the first row of
the current tableau, place $\sigma(i)$ at the end of the first
row. Otherwise, let $m$ be the leftmost element on the first row
that is larger than $\sigma(i)$. Place $\sigma(i)$ in the square
that $m$ occupied, and place $m$ at the end of the second row (in
this case we say that $\sigma(i)$ \emph{bumps} $m$). The
\emph{recording tableau} $Q$ has the same shape as $P$ and is
obtained by placing~$i$ in the position of the
square that was created at step $i$ (when $\sigma(i)$ was
inserted) in the construction of $P$, for all~$i$ from~1 to~$n$.
We write $\mathrm{RSK}(\sigma)=(P,Q)$.

\begin{figure}[hbt]
\epsfig{file=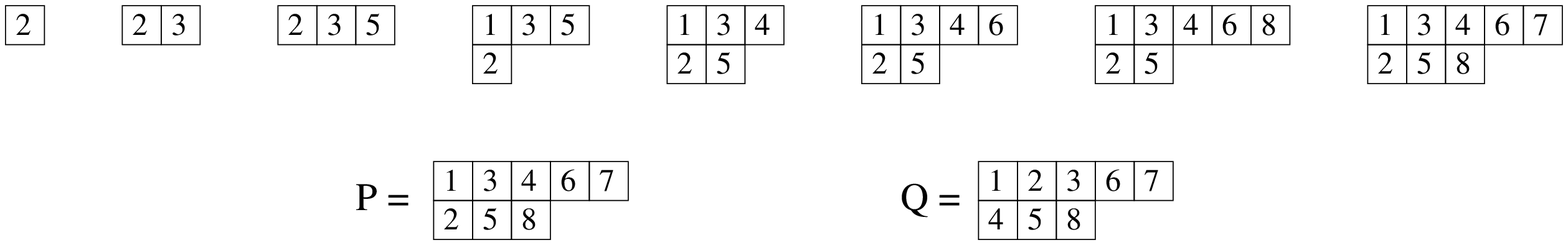,width=6in} \caption{\label{fig:rsk}
Construction of the  RSK correspondence $\mathrm{RSK}(\sigma)=(P,Q)$
for $\sigma=(2,3,5,1,4,6,8,7)$.}
\end{figure}

Now, the first half of the Dyck path $\knu(\sigma)$ is obtained by
adjoining, for $i$ from 1 to $n$, an up-step if $i$ is on the
first row of $P$, and a down-step if it is on the second row. Let
$A$ be the corresponding word of $u$'s and $d$'s. Similarly, let
$B$ be the word obtained from $Q$ in the same way. We define
$\knu(\sigma)$ to be the Dyck path obtained by the concatenation
of the word $A$ and the word $B$ written backwards. For example,
from the tableaux $P$ and $Q$ as in Figure~\ref{fig:rsk} we get the
Dyck path shown in Figure~\ref{fig:ctrt}. The following proposition
summarizes properties of this bijection~$\knu$:

\begin{prop}\label{prop:321} The bijection $\knu:\S_n(321)\longrightarrow\D_n$
satisfies $\fp(\sigma)=\ct(\knu(\sigma))$,
$\exc(\sigma)=\rt(\knu(\sigma))$, and
$\ell(\sigma)=\frac{1}{2}\bigl(n+\he(\knu(\sigma))\bigr)$, for all
$\sigma \in \S_n(321)$.
\end{prop}

Suppose $\mathrm{RSK}(\sigma)=(P,Q)$ for any $\sigma \in S_n$. A
fundamental and highly nontrivial property of the RSK
correspondence is the {\it duality}:
$\mathrm{RSK}(\sigma^{-1})=(Q,P)$ (see e.g.~\cite{EC}). The
classical {\it Schensted's Theorem} states that $\ell(\sigma)$ is
equal to the length of the first row of the tableau~$P$ (and~$Q$).
Both results are used in the proof of Proposition~\ref{prop:321}.

\ms

Let us now define the second bijection~$\krar$, which is essentially
the inverse of the bijection between $\S_n(132)$ and $\D_n$ given by
Krattenthaler \cite{Kra01}, up to reflection of the path from a
vertical line. Following the presentation in Reifegerste
\cite{Rei02}, the path $\krar(\sigma)$ can be described in terms
of the \emph{diagram} of $\sigma$, which is obtained from the
$n\times n$ array representation of $\sigma$ by shading, for each
cross, the cell containing it and the squares that are due south
and due east of it. This gives a bijection between $\S_n(132)$ and
Young diagrams that fit in the shape $(n-1,n-2,\ldots,1)$.
Consider now the path determined by the border of the diagram,
that is, the path with \emph{up} and \emph{right} steps that goes
from the lower-left corner to the upper-right corner of the array,
leaving all the crosses to the right, and staying always as close
to the diagonal connecting these two corners as possible. Define
$\krar(\sigma)$ to be the Dyck path obtained from this path by
reading an up-step every time it goes up and a down-step every
time it goes right.

\begin{figure}[hbt]
\epsfig{file=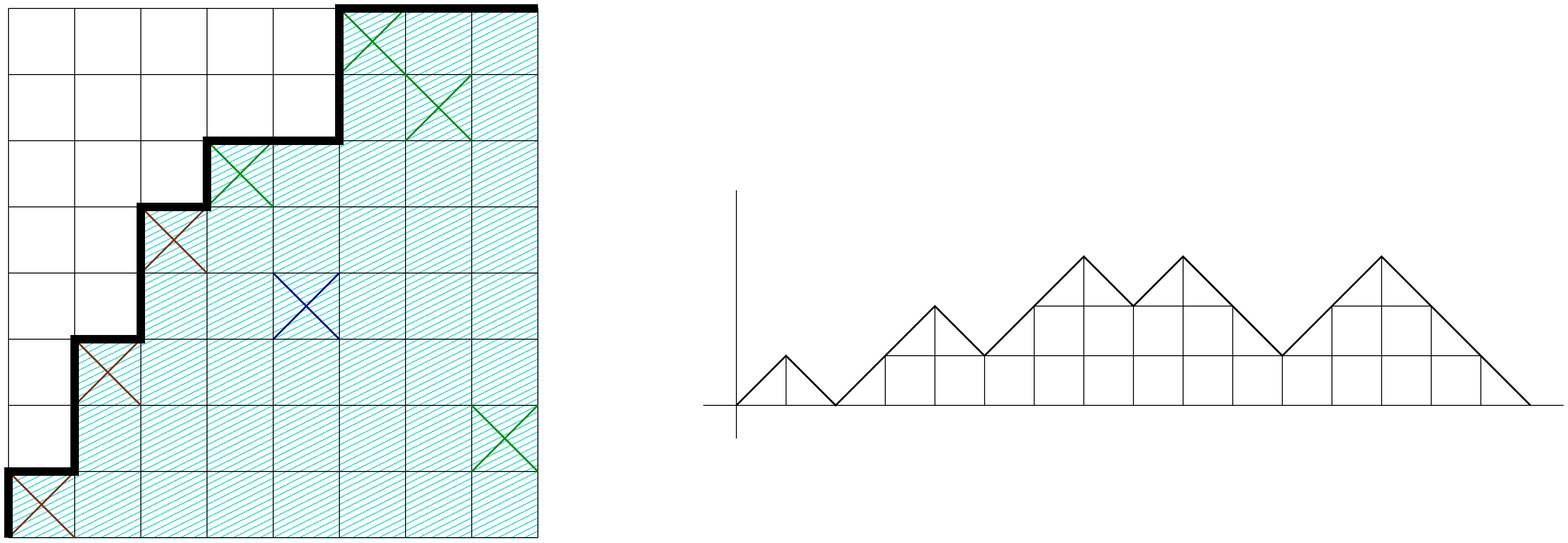,height=3.7cm} \caption{\label{fig:bij2}
The bijection $\krar:(6,7,4,3,5,2,8,1)\mapsto uduuduududduuddd$.}
\end{figure}

\begin{prop}\label{prop:132} The bijection $\krar:\S_n(132)\longrightarrow\D_n$
satisfies that $\fp(\sigma)=\ct(\krar(\sigma))$,
$\exc(\sigma)=\rt(\krar(\sigma))$, and
$\rk(\sigma)=\frac{1}{2}\bigl(n-\he(\krar(\sigma))\bigr)$,
for all $\sigma \in \S_n(132)$.
\end{prop}

\ms

The main result of the paper follows now easily from these two
propositions.

\begin{proof}[Proof of Theorem~\ref{th:new}]
Propositions~\ref{prop:321} and \ref{prop:132} imply that
$\XX=\krar^{-1}\circ\knu$ is a bijection from $\S_n(321)$ to
$\S_n(132)$ which satisfies
$\fp(\XX(\sigma))=\ct(\knu(\sigma))=\fp(\sigma)$,
$\exc(\XX(\sigma))=\rt(\knu(\sigma))=\exc(\sigma)$, and
$$\rk(\XX(\sigma))=\frac{1}{2}\bigl(n-\he(\knu(\sigma))\bigr)
=n-\frac{1}{2}\bigl(n+\he(\knu(\sigma))\bigr)=n-\ell(\sigma)\,.$$
This implies the result.
\end{proof}

\section{Proof of proposition~\ref{prop:132}}\label{sec:pf132}

It can be seen using the diagram representation that $\krar$ maps
fixed points to centered tunnels and excedances to right tunnels.
There is an easy way to recover a permutation $\sigma\in\S_n(321)$
from its diagram: row by row, put a cross in the leftmost shaded
square such that there is exactly one dot in each column.

Instead of looking directly at $\krar(\sigma)$, consider the path
from the lower-left corner to the upper-right corner of the array
of $\sigma$. To each cross we can associate a tunnel in a natural
way. Indeed, each cross produces a decomposition
$\krar(\sigma)=AuBdC$ where $B$ corresponds to the part of the
path above and to the left of the cross. Here $u$ corresponds to
the vertical step directly to the left of the cross, and $d$ to
the horizontal step directly above the cross. According to whether
the cross was to the left of, to the right of, or on the main
diagonal, it will produce respectively a left, right, or centered
tunnel. Thus, fixed points give centered tunnels and excedances give
right tunnels.


Similarly, the rank $\rk(\sigma)$ is the largest $m$
such that an $m\times m$ square fits in the upper-left corner
of the diagram of~ $\sigma$. If we scale the size
of the array so that its diagonal has length $2n$, then the
diagonal of this $\rk(\sigma)\times\rk(\sigma)$ square has length
$2\rk(\sigma)$, and the height of $\krar(\sigma)$ is exactly
$\he(\krar(\sigma))=n-2\rk(\sigma)$.

\section{Proof of proposition~\ref{prop:321}}\label{sec:pf321}

Let us first consider only fixed points in a permutation $\sigma \in \S_n$.
Observe that if $\sigma\in\S_n(321)$ and $\sigma(i)=i$, then
$(\sigma(1),\sigma(2),\ldots,\sigma(i-1))$ is a permutation of
$\{1,2,\ldots,i-1\}$, and
$(\sigma(i+1),\sigma(i+2),\ldots,\sigma(n))$ is a permutation of
$\{i+1,i+2,\ldots,n\}$. Indeed, if $\sigma(j)>i$ for some $j<i$,
then necessarily $\sigma(k)<i$ for some $k>i$, and
$(\sigma(j),\sigma(i),\sigma(k))$ would be an occurrence of 321.

Therefore, when we apply RSK to $\sigma$, the elements
$\sigma(i),\sigma(i+1),\ldots,\sigma(n)$ will never bump any of
the elements $\sigma(1),\sigma(2),\ldots,\sigma(i-1)$. In
particular, the subtableaux of $P$ and $Q$ determined by the
entries that are smaller than $i$ will have both the same shape.
Furthermore,
when the elements greater than $i$ are placed in $P$ and $Q$, the
rows in which they are placed are independent of the subpermutation
$(\sigma(1),\sigma(2),\ldots,\sigma(i-1))$. Note also that
$\sigma(i)$ will never be bumped.

When the Dyck path $\knu(\sigma)$ is built from $P$ and $Q$, this
translates into the fact that the steps corresponding to
$\sigma(i)$ in $P$ and to $i$ in $Q$ will be respectively an
up-step in the first half and a down-step in the second half, both
at the same height and at the same distance from the center of the
path. Besides, the part of the path between them will be itself
the Dyck path corresponding to
$(\sigma(i+1)-i,\sigma(i+2)-i,\ldots,\sigma(n)-i)$. So, the fixed
point $\sigma(i)=i$ determines a centered tunnel in
$\knu(\sigma)$. It is clear that the converse is also true, that
is, every centered tunnel comes from a fixed point. This shows
that $\fp(\sigma)=\ct(\knu(\sigma))$.

\smallskip

Let us now consider excedances in a permutation~$\sigma\in\S_n(321)$.
Our goal is to show that the excedances of $\sigma$ correspond to
right tunnels of $\knu(\sigma)$. The first observation is that we
can assume without loss of generality that $\sigma$ has no fixed
points. Indeed, if $\sigma(i)=i$ is a fixed point of $\sigma$,
then the above reasoning shows that we can decompose
$\knu(\sigma)=AuBdC$, where $AC$ is the Dyck path
$\knu((\sigma(1),\sigma(2),\ldots,\sigma(i-1)))$ and $B$ is a
translation of the Dyck path
$\knu((\sigma(i+1)-i,\ldots,\sigma(n)-i))$. But we have that
$\exc(\sigma)=\exc((\sigma(1),\sigma(2),\ldots,\sigma(i-1)))+\exc((\sigma(i+1)-i,\ldots,\sigma(n)-i))$
and $\rt(AuBdC)=\rt(AC)+\rt(B)$, so in this case the result holds
by induction on the number of fixed points. Note also that the
above argument showed that
$\fp(\sigma)=\fp((\sigma(1),\sigma(2),\ldots,\sigma(i-1)))+\fp((\sigma(i+1)-i,\ldots,\sigma(n)-i))+1$
and $\ct(AuBdC)=\ct(AC)+\ct(B)+1$.

Suppose that $\sigma\in\S_n(321)$ has no fixed points. It is known
that a permutation is 321-avoiding if and only if both the
subsequence determined by its excedances and the one determined by
the remaining elements (in this case, the antiexcedances) are
increasing (see e.g.~\cite{Rei02}). Denote by $X_i:=(i,\sigma(i))$
the crosses of the array representation of $\sigma$.
To simplify the presentation, we will refer indistinctively
to $i$ or $X_i$, hoping
this does not lead to a confusion. For example, we will say
``$X_i$ is an excedance", etc.

Define a matching between excedances and antiexcedances of
$\sigma$ by the following algorithm. Let
$\sigma(i_1)<\sigma(i_2)<\cdots<\sigma(i_k)$ be the excedances of
$\sigma$ and let $\sigma(j_1)<\sigma(j_2)<\cdots<\sigma(j_{n-k})$
be the antiexcedances. \ben
\item Initialize $a:=1$, $b:=1$.
\item Repeat until $a>k$ or $b>n-k$: \begin{itemize}
\item If $i_a>j_b$, then $b:=b+1$. ($X_{j_b}$ is not matched.)
\item Else if $\sigma(i_a)<\sigma(j_b)$, then $a:=a+1$. ($X_{i_a}$ is not
matched.)
\item Else, match $X_{i_a}$ with $X_{j_b}$; $a:=a+1$, $b:=b+1$.
\end{itemize}\een

\begin{figure}[hbt]
\epsfig{file=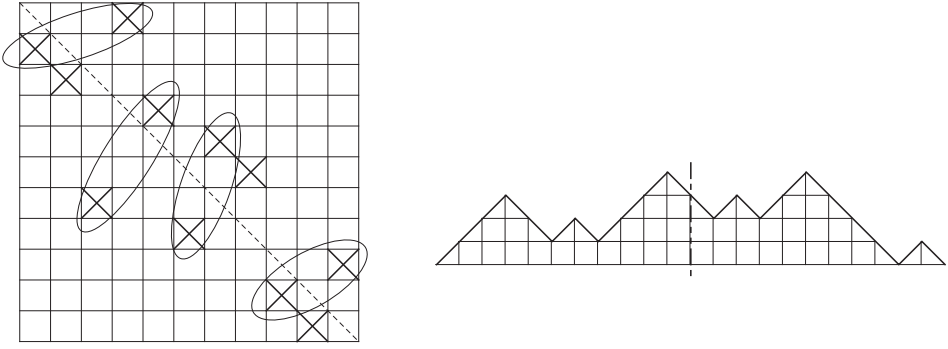,height=2in} \caption{\label{fig:match}
Example of the matching for $\sigma=(4,1,2,5,7,8,3,6,11,9,10)$,
and $\knu(\sigma)$.}
\end{figure}

Now we consider the matched excedances on one hand and the
unmatched ones on the other. We summarize rather technical
results in the following two lemmas, which are proved in
section~\ref{sec:pflemmas}.

\begin{lemma}\label{lemma:matched} The following quantities are equal: \ben
\item the number of matched pairs $(X_i,X_j)$, where $X_i$ is an excedance
and $X_j$ an antiexcedance;
\item the length of the second row of $P$ (or $Q$);
\item the number of right-side tunnels of $\knu(\sigma)$;
\item the number of left-side tunnels of $\knu(\sigma)$;
\item $\frac{1}{2}\bigl(n-\he(\knu(\sigma))\bigr)$;
\item $n-\ell(\sigma)$.
\een\end{lemma}

Note that (5)=(6) implies that
$\ell(\sigma)=\frac{1}{2}\bigl(n+\he(\knu(\sigma))\bigr)$, which is the third
part of Proposition~\ref{prop:321}.

\begin{lemma}\label{lemma:unmatched}
The number of unmatched excedances (resp.
antiexcedances) of $\sigma$ equals the number of right-across
(resp. left-across) tunnels of $\knu(\sigma)$.
\end{lemma}

Since each excedance of $\sigma$ either is part of a matched pair
$(X_i,X_j)$ or is unmatched, lemmas~\ref{lemma:matched} and
\ref{lemma:unmatched} imply that the total number $\exc(\sigma)$
of excedances equals the number of right-side tunnels of
$\knu(\sigma)$ plus the number of right-across tunnels, which is
$\rt(\knu(\sigma))$.

\section{Proofs of the lemmas}\label{sec:pflemmas}

\begin{proof}[Proof of Lemma~\ref{lemma:matched}]
From the descriptions of the RSK algorithm and the matching, it
follows that an excedance $X_i$ and an antiexcedance $X_j$ are
matched with each other precisely when $\sigma(j)$ bumps
$\sigma(i)$ when RSK is performed on $\sigma$, and that these are
the only bumpings that take place. Indeed, an excedance never
bumps anything because it is larger than the elements inserted
before. On the other hand, when an antiexcedance $X_j$ is
inserted, it bumps the smallest element larger than $\sigma(j)$
which has not been bumped yet (which corresponds to an excedance
that has not been matched yet), if such element exists. This
proves the equality (1)=(2).

To see that (2)=(3), observe that right-side tunnels correspond to
up-steps in the right half of $\knu(\sigma)$, which by the
construction of the bijection~$\knu$ correspond to elements on the
second row of $Q$. The equality (3)=(5) follows easily by
counting the number of up-steps and down-steps of the right half
of the path. The equality (4)=(5) is analogous.

Finally, Schensted's Theorem states that the size of the first
row of $P$ equals the length of the longest increasing subsequence
of $\sigma$. This implies that (2)=(6), which completes the proof.
\end{proof}

The reasoning used in the above proof gives a nice
equivalent description of the recording tableau $Q$ in terms of
the array and the matching. Read the rows of the array from top to
bottom. For $i$ from 1 to $n$, place $i$ on the first row of $Q$
if $X_i$ is an excedance or it is unmatched, and place $i$ on the
second row if $X_i$ is a matched antiexcedance. In the
construction of the right half of $\knu(\sigma)$, this translates
into drawing the path from right to left while reading the array
from top to bottom, adjoining an up-step for each matched
antiexcedance and a down-step for each other kind of cross.

To get a similar description of the tableau~$P$, use the fact that the
matching is invariant under transposition (reflection along the
main diagonal) of the array, by the way it was defined.
Recall the duality of the RSK correspondence: if
$\mathrm{RSK}(\sigma)=(P,Q)$, then $\mathrm{RSK}(\sigma^{-1})=(Q,P)$.
Therefore, tableau~$P$ can be obtained by
reading the columns of the array of $\sigma$ from left to right
and, for each column $j$, placing $j$ on the first row of $P$ if
the cross in column $j$ is an antiexcedance or it is unmatched,
and placing $j$ on the second row if the cross is a matched
excedance. Equivalently, the left half of $\knu(\sigma)$, from
left to right, is obtained by reading the array from left to right
and adjoining a down-step for each matched excedance, and an
up-step for each other kind of cross.

In particular, when the left half of the path is constructed in
this way, every matched pair $(X_i,X_j)$ produces an up-step and a
down-step, giving the latter a left-side tunnel. Similarly, in the
construction of the right half of the path, a matched pair gives a
right-side tunnel.

\begin{proof}[Proof of Lemma~\ref{lemma:unmatched}]
It is enough to prove it only for the case of excedances. The case
of antiexcedances follows from it considering $\sigma^{-1}$ and
noticing that $\knu(\sigma^{-1})$ is obtained from $\knu(\sigma)$
by reflecting it into a vertical line.
Let $X_k$ be an unmatched excedance of $\sigma$. We use the above
description of $\knu(\sigma)$ in terms of the array and the
matching. Each cross $X_i$ produces a step $r_i$ in the right half
of the Dyck path and another step $l_i$ in the left half. Crosses
above $X_k$ produce steps to the right of $r_k$, and crosses to
the left of $X_k$ produce steps to the left of $l_k$. In
particular, there are $k-1$ steps to the right of $r_k$, and
$\sigma(k)-1$ steps to the left of $l_k$. Note that since $X_k$ is
an excedance and $\sigma$ is 321-avoiding, all the crosses above
it are also to the left of it.
Consider the crosses that lie to the left of $X_k$.
They can be of the following four kinds:
\smallskip
\begin{itemize}
\item {\it Unmatched excedances $X_i$}.
They will necessarily lie above $X_k$, because the subsequence of
excedances of $\sigma$ is decreasing. Each one of these crosses
contributes an up-step to the left of $l_k$ and down-step to the
right of $r_k$.
\smallskip
\item {\it Unmatched antiexcedances $X_j$}.
They also have to lie above $X_k$, otherwise $X_k$ would be
matched with one of them.
So, each such $X_j$ contributes an up-step to the left of
$l_k$ and down-step to the right of $r_k$.
\smallskip
\item {\it Matched pairs $(X_i,X_j)$} (i.e. $X_i$ is an excedance
and $X_j$ an antiexcedance){\it, where both $X_i$ and $X_j$
lie above $X_k$}.
Both crosses together will contribute an up-step and a
down-step to the left of $l_k$, and an up-step and a down-step
to the right of $r_k$.
\smallskip
\item {\it Matched pairs $(X_i,X_j)$} (i.e. $X_i$ is an excedance and $X_j$
an antiexcedance){\it, where $X_j$ lies below $X_k$}.
The pair will contribute an up-step and a down-step to the left of $l_k$.
However, to the right of $r_k$,
the only contribution will be a down-step produced by $X_i$.
\end{itemize}
\smallskip
Note that there cannot be an antiexcedance $X_j$ to the left of
$X_k$ matched with an excedance to the right of $X_k$, because in
this case $X_j$ would have been matched with $X_k$ by the
algorithm. In the first three cases, the contribution to both
sides of the Dyck path is the same, so that the heights of $r_k$
and $l_k$ are equally affected. But since $\sigma(k)>k$, at least
one of the crosses to the left of $X_k$ must be below it, and this
must be a matched antiexcedance as in the fourth case. This
implies that the step $r_k$ is at a higher $y$-coordinate than
$l_k$. Let $h_k$ be the height of $l_k$. We now show that
$\knu(\sigma)$ has a right-across tunnel at height $h_k$.

Observe that~$h_k$ is the number of unmatched
crosses to the left of $X_k$, and that the height of $r_k$ is the
number of unmatched crosses above $X_k$ (which equals $h_k$) plus
the number of excedances above $X_k$ matched with antiexcedances
below $X_k$. The part of the path between $l_k$ and the middle
always remains at a height greater than $h_k$. This is because the
only possible down-steps in this part can come from matched
excedances $X_i$ to the right of $X_k$, but then such a $X_i$ is
matched with an antiexcedance $X_j$ to the right of $X_k$ but to
the left of $X_i$, which produces an up-step ``compensating" the
down-step associated to $X_i$. Similarly, the part of the path
between $r_k$ and the middle remains at a height greater than
$h_k$. This is because the $h_k$ down-steps to the right of $r_k$
that come from unmatched crosses above $X_k$ don't have a
corresponding up-step in the part of the path between $r_k$ and
the middle. Hence, $l_k$ is the left end of a right-across tunnel,
since the right end of this tunnel is to the right of $r_k$, which
in turn is closer to the right end of $\knu(\sigma)$ than $l_k$ is
to its left end.

It can easily be checked that the converse is also true, namely
that in every right-across tunnel of $\knu(\sigma)$, the step at
its left end corresponds to an unmatched excedance of $\sigma$.
\end{proof}

\section{Final Remarks}\label{sec:remarks}

First, recall the result in~\cite{RSZ02} that the number of
permutations
$\sigma \in \S_n(132)$ (or $\sigma \in \S_n(321)$)
with no fixed points is the {\it Fine number}~$F_n$.
This sequence is most easily defined by its relation
to Catalan numbers:
$$C_n \, = \, 2\,F_n \, + \, F_{n-1} \ \  \text{for} \ n \ge 2, \
\ \text{and} \ F_1=0, \, F_2=1.$$ Although defined awhile ago,
Fine numbers have received much attention in recent years (see a
survey~\cite{DS}). One application of our results are simple
bijections between these two combinatorial interpretations of Fine
numbers and a new one: the set of Dyck paths without centered
tunnels.  In particular, we obtain a bijective proof of the
following result, which follows from~\cite{RSZ02,Eli02}.

\begin{corollary} \label{th:fine}
The number of Dyck paths $D \in \D_n$ without centered
tunnels is equal to~$F_n$.
\end{corollary}

In a different direction, one can extend
Propositions~\ref{prop:321} and~\ref{prop:132} to statistics
$\he_c(D)$ defined as the height at $x=n-c$ of the Dyck path~$D
\in \D_n$, for any $c \in \{0, \pm 1, \pm 2, \dots, \pm (n-1)\}$.
The corresponding statistics in $\S_n(132)$ and in $\S_n(321)$ are
generalizations of the rank of a permutation and the length of the
longest increasing subsequence in a certain subpermutation
of~$\sigma$. The corresponding generalization of
Theorem~\ref{th:new} is straightforward and is left to the reader.

\smallskip

Our final extension has appeared unexpectedly after the results of
this paper have been obtained. We say that a permutation $\sigma
\in \S_n$ is an {\it involution} if~$\sigma = \sigma^{-1}$. In a
recent paper~\cite{DRS} the authors introduce a notion of {\it
refined restricted involutions} by considering ``the number of
fixed points" statistic on involutions avoiding different
patterns $\pi \in \S_3$. They prove the following result:

\begin{theorem}\cite{DRS} \label{th:inv}
\, The number of 321-avoiding involutions $\sigma \in \S_n$ with
$\fp(\sigma)=i$ equals the number of 132-avoiding involutions
$\sigma \in \S_n$ with $\fp(\sigma)=i$, for any $0 \le i \le n$.
\end{theorem}

Let us show that Theorem~\ref{th:inv} follows easily from our
investigation.  Indeed, for every Dyck path $D\in \D_n$ denote by
$D^\ast$ the path obtained by reflection of $D$ from a vertical
line $x=n$.  Now observe that if $\krar(\sigma)=D$, then
$\krar(\sigma^{-1})=D^\ast$.  Similarly, if $\knu(\sigma)=D$, then
$\knu(\sigma^{-1})=D^\ast$ (this follows immediately from the
duality of RSK).  Therefore, $\sigma\in \S_n(321)$ is an
involution if and only if so is $\Theta(\sigma) \in \S_n(132)$,
which implies the result. Furthermore, we obtain the following
extension of Theorem~\ref{th:inv}:

\begin{theorem}\label{th:new-inv}
The number of 321-avoiding involutions $\sigma \in \S_n$
with $\fp(\sigma)=i$, $\exc(\sigma)=j$ and $\ell(\sigma)=k$ equals
the number of 132-avoiding involutions $\sigma \in \S_n$
with $\fp(\sigma)=i$, $\exc(\sigma)=j$ and $\rk(\sigma)=n-k$, for
any $0 \le i,j,k\le n$.
\end{theorem}

\smallskip

Finally, a few questions and open problems.  First, it would be
nice to obtain a ``more philosophical'' proof of
Theorems~\ref{th:old} and~\ref{th:new}, to see if this is more
than a (proven) coincidence.  Is there any reason to believe that
refined restricted permutations are equinumerous for some larger
patterns?  What happens, for example, for patterns $\pi \in \S_4$?

Second, note that the RSK correspondence arose naturally in our
investigation as well as in~\cite{DRS}.  Is there a more general
result on pattern avoidance which uses RSK to a larger extend?  Is
there a general result in representation theory of $\S_n$ which
might explain Theorem~\ref{th:new}?  We hope the reader is as
puzzled as we are at this point.


\subsection*{Acknowledgements}

We would like to thank Richard Stanley for suggesting
the problem of enumerating excedances in pattern-avoiding
permutations and for helpful conversations.
The first author was partially supported by the MAE.
The second author was supported by the NSA and the NSF.

\ms

\end{document}